\documentclass[conference]{IEEEtran}
\IEEEoverridecommandlockouts

\usepackage{amsmath,amssymb}
\usepackage{graphicx}
\usepackage{cite}
\usepackage{bm}
\usepackage{url}
\usepackage{siunitx}
\usepackage{color}
\usepackage[normalem]{ulem}
\usepackage{tabularx}
\usepackage{booktabs}
\usepackage{float}
\usepackage{comment}
\usepackage{tikz}
\usepackage{lettrine}

\title{A Mode-Matching Approach to the Design\\
of RIS-Aided Communications}
\author{
\IEEEauthorblockN{
Ahmed Najjar, 
Hajar El~Hassani,  
Marco Di~Renzo, \IEEEmembership{Fellow,~IEEE}, 
Kezhi~Wang, \IEEEmembership{Senior~Member,~IEEE}, \\
Merouane Debbah, \IEEEmembership{Fellow,~IEEE}
}
\thanks{This work was supported in part by the European Union through the Horizon Europe project COVER under grant agreement number 101086228, the Horizon Europe project UNITE under grant agreement number 101129618, the Horizon Europe project INSTINCT under grant agreement number 101139161, and the Horizon Europe project TWIN6G under grant agreement number 101182794, as well as by the Agence Nationale de la Recherche (ANR) through the France 2030 project ANR-PEPR Networks of the Future under grant agreements NF-Founds 22-PEFT-0010 and NF-YACARI 22-PEFT-0005, and by the CHIST-ERA project PASSIONATE under grant agreements CHIST-ERA-22-WAI-04 and ANR-23-CHR4-0003-01.}

\thanks{A. Najjar, is  with Universit\'e Paris-Saclay, CNRS, CentraleSup\'elec, Laboratoire des Signaux et Syst\`emes, 3 Rue Joliot-Curie, 91192 Gif-sur-Yvette, France (ahmed.najjar@centralesupelec.fr).} 
\thanks{M. Di Renzo is with Universit\'e Paris-Saclay, CNRS, CentraleSup\'elec, Laboratoire des Signaux et Syst\`emes, 3 Rue Joliot-Curie, 91192 Gif-sur-Yvette, France. (marco.di-renzo@universite-paris-saclay.fr), and with King's College London, Department of Engineering - Centre for Telecommunications Research, WC2R 2LS London, United Kingdom (marco.di\_renzo@kcl.ac.uk).}
\thanks{H. El Hassani is with ETIS, UMR 8051,
CY Cergy Paris Universit\'e, ENSEA, CNRS, F-95000 Cergy, France. (hajar.el-hassani@ensea.fr).}

\thanks{K. Wang is with the Department of Computer Science, Brunel University of London, United Kingdom (kezhi.wang@brunel.ac.uk).} 
\thanks{M. Debbah is with KU 6G Research Center, Khalifa University of Science and Technology, P O Box 127788, Abu Dhabi, UAE (merouane.debbah@ku.ac.ae)}

}

\begin{document}

\maketitle

\begin{abstract}
Reconfigurable intelligent surface (RIS) is an emerging technology for application to wireless communications. In this paper, we consider the problem of anomalous reflection and model the RIS as a periodic surface impedance boundary. We utilize the mode-matching method and Floquet’s expansion representation to compute the field reflected from a spatially periodic RIS, and evaluate the performance versus implementation complexity tradeoffs of RIS-aided communications based on the global 
design criterion. This allows us to maximize the power reflected towards the intended direction of propagation, while minimizing the power reradiated towards undesired directions of propagation. In addition, we discuss the advantages of the proposed electromagnetically consistent approach to the design of RIS-aided wireless systems.
\end{abstract}
\section{Introduction}
Reconfigurable intelligent surface (RIS) has recently emerged as a promising technology for shaping wireless propagation environments in a nearly passive way. By carefully engineering programmable meta‑atoms whose electromagnetic (EM) response can be tuned in real time, RIS enables the realization of smart radio environments capable of manipulating incident waves in a controlled fashion \cite{di2020smart},\cite{di2022communication}. A fundamental functionality of RIS is anomalous reflection, where an incoming plane wave is redirected towards a specified non‑specular direction. 

Most RIS‑aided communication models rely on locally periodic phase‑gradient metasurfaces, often derived from the generalized Snell’s law or reflectarray theory \cite{liu2023reflectarrays}. These local design approaches assign a spatially varying phase shift to each RIS element and implicitly assume each of them is designed by assuming locally periodic boundary conditions \cite{di2022communication}. While attractive due to their simplicity, these methods are known to suffer from significant efficiency degradation when steering waves towards large reflection angles, mainly due to the unavoidable excitation of parasitic Floquet harmonics \cite{mohammadi2016wave,asadchy2016perfect,diaz2017generalized,epstein2016synthesis,kwon2018lossless,asadchy2017eliminating,he2022perfect}. These effects arise because local phase designs neglect the electromagnetic coupling between adjacent elements, ignore the contribution of evanescent modes, and do not enforce global power conservation.

To overcome these limitations, several works have introduced electromagnetically consistent and global design frameworks for RIS based on inhomogeneous surface‑impedance boundaries \cite{di2022communication,diaz2021macroscopic}. These models explicitly account for the interaction between all propagating and evanescent Floquet harmonics and reveal the fundamental differences between local and global surface designs in terms of power flow, scattering behavior, and achievable efficiency. More recently, scalable EM‑consistent synthesis methods \cite{vuyyuru2024modeling} and global impedance‑optimization algorithms \cite{shabir2025electromagnetically} and \cite{kosulnikov2024experimental} have further emphasized the importance of incorporating Maxwell’s equations and physical constraints into RIS design. Within this context, modeling RIS as surface‑impedance boundaries offers a physically-consistent framework for analyzing anomalous reflection. The mode‑matching technique provides an exact decomposition of the reflected field into Floquet harmonics and enables the computation of their amplitudes for any prescribed periodic surface‑impedance profile, which provides a rigorous characterization of both propagating and evanescent waves \cite{hwang2012periodic}. 

In this work, we analyze anomalous reflection from a periodic RIS using an electromagnetically consistent formulation. We derive the surface‑impedance expression from the tangential fields via a Floquet expansion and compute the reflected Floquet coefficients using the mode-matching technique. We then evaluate three representative impedance profiles to quantify the performance versus implementation complexity tradeoffs. The results highlight the importance of physically consistent impedance designs for efficient and directive anomalous reflection in RIS‑aided communication systems. 
\section{System Model}
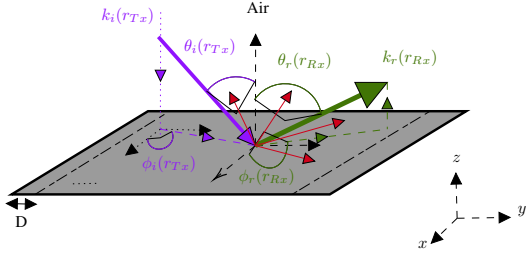
\begin{figure}[t]
    \centering
    \resizebox{0.8\columnwidth}{!}{%
        \tikzset{every picture/.style={line width=0.55pt}}
        
        \begin{tikzpicture}[x=0.75pt,y=0.75pt,yscale=-1,xscale=1]
        
        \draw  [fill={rgb, 255:red, 155; green, 155; blue, 155 }  ,fill opacity=1 ][line width=1.5]  (251,103) -- (509.3,103) -- (398.6,170.6) -- (140.3,170.6) -- cycle ;
        \draw  [dash pattern={on 4.5pt off 4.5pt}]  (338.2,131) -- (387.4,130.62) ;
        \draw [shift={(390.4,130.6)}, rotate = 179.56] [fill={rgb, 255:red, 0; green, 0; blue, 0 }  ][line width=0.08]  [draw opacity=0] (8.93,-4.29) -- (0,0) -- (8.93,4.29) -- cycle   ;
        \draw  [dash pattern={on 4.5pt off 4.5pt}]  (338.2,131) -- (305.95,157.34) ;
        \draw [shift={(304.4,158.6)}, rotate = 320.77] [color={rgb, 255:red, 0; green, 0; blue, 0 }  ][line width=0.75]   (10.93,-3.29) .. controls (6.95,-1.4) and (3.31,-0.3) .. (0,0) .. controls (3.31,0.3) and (6.95,1.4) .. (10.93,3.29)  ;
        \draw  [dash pattern={on 4.5pt off 4.5pt}]  (338.2,131) -- (338.39,43.6) ;
        \draw [shift={(338.4,40.6)}, rotate = 90.13] [fill={rgb, 255:red, 0; green, 0; blue, 0 }  ][line width=0.08]  [draw opacity=0] (8.93,-4.29) -- (0,0) -- (8.93,4.29) -- cycle   ;
        \draw [color={rgb, 255:red, 65; green, 117; blue, 5 }  ,draw opacity=1 ][line width=3]   (338.2,131) -- (439.89,82.2) ;
        \draw [shift={(445.3,79.6)}, rotate = 154.36] [fill={rgb, 255:red, 65; green, 117; blue, 5 }  ,fill opacity=1 ][line width=0.08]  [draw opacity=0] (23.75,-11.41) -- (0,0) -- (23.75,11.41) -- cycle   ;
        \draw [color={rgb, 255:red, 144; green, 19; blue, 254 }  ,draw opacity=1 ][line width=2.25]   (259.4,42.6) -- (334.87,127.27) ;
        \draw [shift={(338.2,131)}, rotate = 228.29] [fill={rgb, 255:red, 144; green, 19; blue, 254 }  ,fill opacity=1 ][line width=0.08]  [draw opacity=0] (14.29,-6.86) -- (0,0) -- (14.29,6.86) -- cycle   ;
        \draw  [dash pattern={on 0.84pt off 2.51pt}]  (261,117.6) -- (299,118.53) ;
        \draw [shift={(302,118.6)}, rotate = 181.4] [fill={rgb, 255:red, 0; green, 0; blue, 0 }  ][line width=0.08]  [draw opacity=0] (8.93,-4.29) -- (0,0) -- (8.93,4.29) -- cycle   ;
        \draw  [dash pattern={on 0.84pt off 2.51pt}]  (261,117.6) -- (235.4,136.8) ;
        \draw [shift={(233,138.6)}, rotate = 323.13] [fill={rgb, 255:red, 0; green, 0; blue, 0 }  ][line width=0.08]  [draw opacity=0] (8.93,-4.29) -- (0,0) -- (8.93,4.29) -- cycle   ;
        \draw [color={rgb, 255:red, 144; green, 19; blue, 254 }  ,draw opacity=1 ] [dash pattern={on 0.84pt off 2.51pt}]  (261.3,22.6) -- (261,124.6) ;
        \draw [shift={(261.14,78.6)}, rotate = 270.17] [fill={rgb, 255:red, 144; green, 19; blue, 254 }  ,fill opacity=1 ][line width=0.08]  [draw opacity=0] (8.93,-4.29) -- (0,0) -- (8.93,4.29) -- cycle   ;
        \draw [color={rgb, 255:red, 144; green, 19; blue, 254 }  ,draw opacity=1 ] [dash pattern={on 4.5pt off 4.5pt}]  (261,117.6) -- (338.2,131) ;
        \draw [shift={(304.53,125.16)}, rotate = 189.85] [fill={rgb, 255:red, 144; green, 19; blue, 254 }  ,fill opacity=1 ][line width=0.08]  [draw opacity=0] (8.93,-4.29) -- (0,0) -- (8.93,4.29) -- cycle   ;
        \draw  [draw opacity=0] (270.24,119.63) .. controls (270.36,119.84) and (270.48,120.07) .. (270.58,120.3) .. controls (272.62,124.83) and (269.98,130.44) .. (264.69,132.81) .. controls (259.4,135.19) and (253.46,133.44) .. (251.42,128.9) .. controls (250.92,127.8) and (250.71,126.64) .. (250.73,125.48) -- (261,124.6) -- cycle ; 
        \draw  [color={rgb, 255:red, 144; green, 19; blue, 254 }  ,draw opacity=1 ] (270.24,119.63) .. controls (270.36,119.84) and (270.48,120.07) .. (270.58,120.3) .. controls (272.62,124.83) and (269.98,130.44) .. (264.69,132.81) .. controls (259.4,135.19) and (253.46,133.44) .. (251.42,128.9) .. controls (250.92,127.8) and (250.71,126.64) .. (250.73,125.48) ;  
        \draw [color={rgb, 255:red, 65; green, 117; blue, 5 }  ,draw opacity=1 ] [dash pattern={on 4.5pt off 4.5pt}]  (445.3,79.6) -- (445.3,117.6) ;
        \draw [shift={(445.3,92.1)}, rotate = 90] [fill={rgb, 255:red, 65; green, 117; blue, 5 }  ,fill opacity=1 ][line width=0.08]  [draw opacity=0] (8.93,-4.29) -- (0,0) -- (8.93,4.29) -- cycle   ;
        \draw [color={rgb, 255:red, 65; green, 117; blue, 5 }  ,draw opacity=1 ] [dash pattern={on 4.5pt off 4.5pt}]  (338.2,131) -- (445.3,117.6) ;
        \draw [shift={(396.71,123.68)}, rotate = 172.87] [fill={rgb, 255:red, 65; green, 117; blue, 5 }  ,fill opacity=1 ][line width=0.08]  [draw opacity=0] (8.93,-4.29) -- (0,0) -- (8.93,4.29) -- cycle   ;
        \draw  [draw opacity=0] (298.8,86.8) .. controls (301.52,83.25) and (305.08,80.25) .. (309.36,78.14) .. controls (318.69,73.55) and (329.23,74.24) .. (337.61,79.06) -- (322.62,105.05) -- cycle ; 
        \draw  [color={rgb, 255:red, 144; green, 19; blue, 254 }  ,draw opacity=1 ] (298.8,86.8) .. controls (301.52,83.25) and (305.08,80.25) .. (309.36,78.14) .. controls (318.69,73.55) and (329.23,74.24) .. (337.61,79.06) ;  
        \draw  [draw opacity=0] (337.65,93.43) .. controls (345.29,82.51) and (359.59,77.67) .. (372.67,82.52) .. controls (382.97,86.34) and (389.92,95.22) .. (391.75,105.3) -- (362.23,110.64) -- cycle ; 
        \draw  [color={rgb, 255:red, 65; green, 117; blue, 5 }  ,draw opacity=1 ] (337.65,93.43) .. controls (345.29,82.51) and (359.59,77.67) .. (372.67,82.52) .. controls (382.97,86.34) and (389.92,95.22) .. (391.75,105.3) ;  
        \draw  [draw opacity=0] (364.48,128.11) .. controls (364.12,138.66) and (359.71,147.01) .. (352.58,148.83) .. controls (345.9,150.54) and (338.54,146.05) .. (333.22,137.99) -- (345.15,119.76) -- cycle ; 
        \draw  [color={rgb, 255:red, 65; green, 117; blue, 5 }  ,draw opacity=1 ] (364.48,128.11) .. controls (364.12,138.66) and (359.71,147.01) .. (352.58,148.83) .. controls (345.9,150.54) and (338.54,146.05) .. (333.22,137.99) ;  
        \draw  [dash pattern={on 3.75pt off 3pt on 7.5pt off 1.5pt}]  (480,103) -- (366.3,171.6) ;
        \draw  [dash pattern={on 3.75pt off 3pt on 7.5pt off 1.5pt}]  (264.3,105.6) -- (159.3,170.6) ;
        \draw  [dash pattern={on 0.84pt off 2.51pt}]  (190,162) -- (210.3,161.6) ;
        \draw    (139,178.07) -- (157.3,178.53) ;
        \draw [shift={(160.3,178.6)}, rotate = 181.41] [fill={rgb, 255:red, 0; green, 0; blue, 0 }  ][line width=0.08]  [draw opacity=0] (8.93,-4.29) -- (0,0) -- (8.93,4.29) -- cycle   ;
        \draw [shift={(136,178)}, rotate = 1.41] [fill={rgb, 255:red, 0; green, 0; blue, 0 }  ][line width=0.08]  [draw opacity=0] (8.93,-4.29) -- (0,0) -- (8.93,4.29) -- cycle   ;
        \draw  [dash pattern={on 4.5pt off 4.5pt}]  (502,190) -- (542,189.35) ;
        \draw [shift={(545,189.3)}, rotate = 179.07] [fill={rgb, 255:red, 0; green, 0; blue, 0 }  ][line width=0.08]  [draw opacity=0] (8.93,-4.29) -- (0,0) -- (8.93,4.29) -- cycle   ;
        \draw  [dash pattern={on 4.5pt off 4.5pt}]  (502,190) -- (501.08,156.3) ;
        \draw [shift={(501,153.3)}, rotate = 88.44] [fill={rgb, 255:red, 0; green, 0; blue, 0 }  ][line width=0.08]  [draw opacity=0] (8.93,-4.29) -- (0,0) -- (8.93,4.29) -- cycle   ;
        \draw  [dash pattern={on 4.5pt off 4.5pt}]  (502,190) -- (483.16,208.21) ;
        \draw [shift={(481,210.3)}, rotate = 315.97] [fill={rgb, 255:red, 0; green, 0; blue, 0 }  ][line width=0.08]  [draw opacity=0] (8.93,-4.29) -- (0,0) -- (8.93,4.29) -- cycle   ;
        \draw [color={rgb, 255:red, 208; green, 2; blue, 27 }  ,draw opacity=1 ]  (338.2,131) -- (365.75,89.11) ;
        \draw [shift={(367.4,86.6)}, rotate = 123.33] [fill={rgb, 255:red, 208; green, 2; blue, 27 }  ,fill opacity=1 ][line width=0.08]  [draw opacity=0] (8.93,-4.29) -- (0,0) -- (8.93,4.29) -- cycle   ;
        \draw [color={rgb, 255:red, 208; green, 2; blue, 27 }  ,draw opacity=1 ]  (338.2,131) -- (407.51,111.42) ;
        \draw [shift={(410.4,110.6)}, rotate = 164.22] [fill={rgb, 255:red, 208; green, 2; blue, 27 }  ,fill opacity=1 ][line width=0.08]  [draw opacity=0] (8.93,-4.29) -- (0,0) -- (8.93,4.29) -- cycle   ;
        \draw [color={rgb, 255:red, 208; green, 2; blue, 27 }  ,draw opacity=1 ]  (338.2,131) -- (383.5,142.84) ;
        \draw [shift={(386.4,143.6)}, rotate = 194.65] [fill={rgb, 255:red, 208; green, 2; blue, 27 }  ,fill opacity=1 ][line width=0.08]  [draw opacity=0] (8.93,-4.29) -- (0,0) -- (8.93,4.29) -- cycle   ;
        \draw [color={rgb, 255:red, 208; green, 2; blue, 27 }  ,draw opacity=1 ]  (338.2,131) -- (318.7,90.31) ;
        \draw [shift={(317.4,87.6)}, rotate = 64.39] [fill={rgb, 255:red, 208; green, 2; blue, 27 }  ,fill opacity=1 ][line width=0.08]  [draw opacity=0] (8.93,-4.29) -- (0,0) -- (8.93,4.29) -- cycle   ;
        
        \draw (142,187) node [anchor=north west][inner sep=0.75pt]   [align=left] {D}; 
        \draw (469,207) node [anchor=north west][inner sep=0.75pt]   [align=left] {$x$};
        \draw (550,179) node [anchor=north west][inner sep=0.75pt]   [align=left] {$y$};
        \draw (497,137) node [anchor=north west][inner sep=0.75pt]   [align=left] {$z$};
        \draw (279,42.4) node [anchor=north west][inner sep=0.75pt][color={rgb, 255:red, 144; green, 19; blue, 254 }]  {$\theta _{i}( r_{T}{}_{x})$};
        \draw (355,54.4) node [anchor=north west][inner sep=0.75pt][color={rgb, 255:red, 65; green, 117; blue, 5 }]  {$\theta _{r}( r_{R}{}_{x})$};
        \draw (245,137.4) node [anchor=north west][inner sep=0.75pt][color={rgb, 255:red, 144; green, 19; blue, 254 }]  {$\phi _{i}( r_{T}{}_{x})$};
        \draw (323.3,148.2) node [anchor=north west][inner sep=0.75pt][color={rgb, 255:red, 65; green, 117; blue, 5 }]  {$\phi _{r}( r_{R}{}_{x})$};
        \draw (212,17.4) node [anchor=north west][inner sep=0.75pt][color={rgb, 255:red, 144; green, 19; blue, 254 }]  {$k_{i}( r_{T}{}_{x})$};
        \draw (440,51.4) node [anchor=north west][inner sep=0.75pt][color={rgb, 255:red, 65; green, 117; blue, 5 }]  {$k_{r}( r_{R}{}_{x})$};
        \draw (330,13) node [anchor=north west][inner sep=0.75pt]   {Air};
        
        \end{tikzpicture}
    }%
    \caption{Geometrical configuration of the RIS. The surface lies in the $xy$-plane, is periodic along $y$, and is modeled as a surface-impedance boundary illuminated by a transverse electric (TE)-polarized plane wave.}
    \label{fig:system_model}
\end{figure}

We consider the system depicted in Figure.~\ref{fig:system_model}, which consists of a transmitter (Tx) located at the point $\mathbf{r_{Tx}}$, a receiver (Rx) located at the point $\mathbf{r_{Rx}}$, and an RIS, defined in a three‑dimensional Cartesian coordinate system $(x,y,z)$. The RIS is modeled as a planar surface of negligible thickness compared to the wavelength $\lambda$, and it is placed in the $xy$-plane, i.e., $z=0$, with its center located at the origin.  
The RIS acts as a semi-passive, programmable surface that reshapes the propagation environment. The transmitter is assumed to be located in the far-field of both
the RIS and the receiver, so the signal impinging on the RIS can be
approximated as a plane wave. The incident wave arrives, in the region $z>0$, from an elevation angle $\theta_i$ and azimuth angle $\phi_i$, and is reflected towards the receiver, located in the direction ($\theta_r,\phi_r$). The wave propagates at the speed of light $c$ in free space characterized by permittivity $\varepsilon_0$, permeability $\mu_0$, and free-space impedance $\eta_0=\sqrt{\mu_0/\varepsilon_0}$. To simplify the analysis, the propagation is restricted to the $yz$-plane by assuming $\phi_i = \phi_r = \pi/2$.

Under this assumption, the RIS can be modeled as a one‑dimensional (1D) periodic structure, which extends infinitely along the $x$-axis and exhibits periodicity along the $y$-axis. The spatial period required to redirect the incident wave from $\theta_i$ to $\theta_r$ is 
\begin{equation}
D = \frac{\lambda}{\left|\,\sin \theta_r-\sin \theta_i\,\right|}
\end{equation}

Since the variation is along the $y$-direction, the RIS can be fully characterized by a periodic spatially dependent surface impedance $Z_s(y)$, which represents the programmable electromagnetic response implemented by the RIS elements. Note that $D$ is defined for specific incidence and reflection design angles ($\theta_{id},\theta_{rd}$) that the RIS is engineered to support. In this work we consider $\theta_{id}=\theta_i$ and $\theta_{rd}=\theta_r$.

Let $\mathbf{k}_i$ and $\mathbf{k}_r$ be the wavevectors of the incident and reflected plane waves. These vectors specify, respectively, the direction from which the transmitter illuminates the RIS and the direction towards which the RIS is engineered to reradiate the maximum power. Both wavevectors have  magnitude equal to the free‑space wavenumber $k= \omega \sqrt{\mu_0 \varepsilon_0}$, where $\omega$ is the angular frequency. The wavevectors are defined as follows \cite{di2022communication}:
\begin{align}
\mathbf{k}_i & = k\left( \sin\theta_i\,\hat{\mathbf{y}} -\cos\theta_i\,\hat{\mathbf{z}} \right) \\
\mathbf{k}_r&= k\left( \sin\theta_r\,\hat{\mathbf{y}} +  \cos\theta_r\,\hat{\mathbf{z}} \right)
\end{align}

In what follows, the time dependency $e^{j\omega t}$, where $j$ is the imaginary unit, is assumed and omitted for brevity. Hence, considering transverse electric (TE) polarization (i.e., the electric field is $x$-polarized), the incident electric field and the corresponding magnetic field obtained from Maxwell’s
equations can be written as 
\begin{align}
\mathbf{E}^{\mathrm{inc}}(y,z)
&=
A_0 e^{-j(k_y y - k_z z)} \hat{\mathbf{x}} \\
\mathbf{H}^{\mathrm{inc}}(y,z)
&=
-\frac{A_0}{\eta_0}
e^{-j(k_y y - k_z z)} \left( \cos\theta_i \ \hat{\mathbf{y}}
+
\sin\theta_i \ \hat{\mathbf{z}}
\right)
\end{align}
where $A_0$ is the complex amplitude of the incident electric field, $k_y = k \sin\theta_i$ and $\
k_z = k \cos\theta_i$. Since the RIS is modeled as an electrically thin surface lying in the $xy$-plane, only the tangential components of the incident fields along $\hat{\mathbf{x}}$ and $\hat{\mathbf{y}}$ are relevant. The electric field is tangential along $\hat{\mathbf{x}}$, whereas the magnetic field has both tangential and normal components. Extracting the tangential components gives
\begin{align}
\mathbf{E}_t^{\mathrm{inc}}(y,z)
&=
A_0 \ e^{-j(k_y y - k_z z)} \ \hat{\mathbf{x}} \\
\mathbf{H}_t^{\mathrm{inc}}(y,z)
&=-\frac{\cos\theta_i}{\eta_0} \ A_0 \
e^{-j(k_y y - k_z z)} \ \hat{\mathbf{y}}
\end{align}

We model the electromagnetic response of the RIS using a surface‑impedance representation.
Because of the periodicity of the RIS along the $y$-direction, the reflected field satisfies Floquet’s theorem \cite{hwang2012periodic}, where the tangential reflected electric field can be expressed as a superposition of plane waves (Floquet harmonics) determined by the tangential wavenumber
\begin{equation}
k_{y,n} = k_y + \frac{2\pi}{D} n, \ \ \ n \in \mathbb{Z} \label{kxn}
\end{equation}
and the propagation direction
\begin{equation}
\theta_{\mathrm{r},n} = \arcsin\!\left(\frac{k_{y,n}}{k}\right), \ n \in \mathbb{Z}
\end{equation}
Each harmonic is associated with a longitudinal wavenumber $k_{z,n}$ that depends on the propagation mode
\begin{equation}
k_{z,n} =
\begin{cases}
\sqrt{k^2 - k_{y,n}^2}, & |k_{y,n}|\le k \ \ \ \  \text{propagating mode}\\[1mm]
j\sqrt{k_{y,n}^2 - k^2} & |k_{y,n}|> k \ \ \ \ \text{evanescent mode}
\end{cases}
\end{equation}

The tangential reflected electric and magnetic fields can be written as a superposition of Floquet harmonics as
\begin{align}
\mathbf{E}_t^{\mathrm{ref}}(y,z)
&= A_0
\sum_{n=-\infty}^{+\infty}
B_n\, 
e^{-j\left(k_{y,n} y + k_{z,n} z\right)} \ \hat{\mathbf{x}}
\label{eq:Escat_TM} \\
\mathbf{H}_t^{\mathrm{ref}}(y,z)
&= A_0
\sum_{n=-\infty}^{+\infty}
Y_n B_n\, 
e^{-j\left(k_{y,n} y + k_{z,n} z\right)} \ \hat{\mathbf{y}} 
\end{align}
where $B_n$ denotes the complex normalized amplitude (reflection coefficient) of the $n$th
reflected Floquet harmonic, and $Y_n = \frac{k_{z,n}}{\eta_0 k}$ is the associated TE wave admittance.

To gain insight for RIS design, we characterize it through the surface impedance $Z_s(y)$, which describes the electromagnetic response of the surface. Since the RIS is modeled as an electrically thin impedance boundary placed at $z=0$, its interaction with the incident wave is determined solely by the total scalar tangential electric and magnetic fields evaluated just above the surface, i.e. $z=0^+$, obtained by summing the incident and reflected contributions as
\begin{align}
E_t^{\mathrm{tot}}(y,0^+)
&=
A_0 \left[ e^{-j k_y y}
+
\sum_{n=-\infty}^{+\infty}
B_n e^{-j k_{y,n} y} \right]
\label{eq:Ex_total_TM}\\
H_t^{\mathrm{tot}}(y,0^+)
& =\! A_0 \!\left[-
\frac{\cos \theta_i}{\eta_0}\, e^{-jk_y y}+\sum_{n=-\infty}^{+\infty}
Y_n B_n\, e^{-jk_{y,n} y} \right]
\label{eq:Hy_total_TM}
\end{align}

The relationship between these total scalar tangential fields is enforced by the surface‑impedance boundary condition as follows:

\begin{equation}
E_{t}^{\mathrm{tot}}(y,0^+)
= - \,Z_{s}(y)\, H_{t}^{\mathrm{tot}}(y,0^+) \label{eq:Zs_def_TM}
\end{equation}
Substituting the total scalar-field expressions \eqref{eq:Ex_total_TM}–\eqref{eq:Hy_total_TM} in \eqref{eq:Zs_def_TM} yields the expression of surface impedance as
\begin{equation}
Z_s(y)
= \,
\frac{1
+ \displaystyle\sum_{n=-\infty}^{+\infty} B_n \Phi_n(y)
}{
\cos\theta_i/\eta_0
- \displaystyle\sum_{n=-\infty}^{+\infty} Y_n B_n \Phi_n(y)
}\label{eq:Zs_final_TM}
\end{equation}
where $
\Phi_n(y) = e^{-j\left(k_{y,n}-k_y\right)y}$. The surface-impedance in \eqref{eq:Zs_final_TM} depends on the unknown Floquet coefficients $B_n$.  These coefficients are determined, for a given surface impedance profile, using the mode-matching technique, which quantifies the resulting power distribution among the reflected Floquet harmonics.

\section{Mode-Matching Based Electromagnetic Model}

To determine the reflected Floquet coefficients generated by a periodic RIS, we adopt the mode-matching technique\cite{hwang2012periodic}. This approach begins by expanding both the surface impedance and the tangential fields in terms of their Floquet series.
Since the surface impedance $Z_s(y)$ is periodic, with period $D$, it admits the Fourier expansion
\begin{equation}
Z_s(y)
=
\sum_{p=-\infty}^{+\infty}z_p\, e^{-j\frac{2\pi}{D} p \, y}
\label{eq:Zs_series}
\end{equation}
where $z_p$ is the $p$th Fourier coefficient. In practice, the infinite Floquet series in \eqref{eq:Ex_total_TM}-\eqref{eq:Hy_total_TM} and \eqref{eq:Zs_series} are truncated to a finite number of harmonics, since only the propagating modes and a limited number of slowly decaying evanescent
modes have a significant contribution to the field. Let $N$ denote the number of harmonics retained on each side of the
fundamental order at $n = 0$. Following the transmission-line analogy for periodically loaded surfaces \cite{hwang2012periodic},
we obtain
\begin{equation}
\mathbf{b}
=
\boldsymbol{\Gamma}\, \mathbf{a}
\label{eq:bGamma}
\end{equation}
where the vectors $\mathbf{a}$ and
$\mathbf{b}$ contain the incident and reflected modal amplitudes, respectively, defined as
\begin{align}
\mathbf{a} &= \bigl[\, 0,\; \ldots,\; 0,\; A_0,\;0,\; \ldots,\; 0 \,\bigr]^{T} \\[4pt]
\mathbf{b} &= A_0 \bigl[\, B_{-N},\; \ldots,\;B_0, \; \ldots,\; B_{N} \,\bigr]^{T}
\end{align}
and $\boldsymbol{\Gamma}$ is the reflection matrix for the truncated set of Floquet harmonics given by \cite{hwang2012periodic}
\begin{equation}
\boldsymbol{\Gamma}
=
\left( \mathbf{I} + \mathbf{Z}_{s}\mathbf{Y}_{a} \right)^{-1}
\left( \mathbf{Z}_{s}\mathbf{Y}_{a} - \mathbf{I} \right)
\label{eq:Gamma_final}
\end{equation}
where $\mathbf{I}$ is the identity matrix, $\mathbf{Y}_a\in\mathbb{C}^{(2N+1)\times (2N+1)}$ is the diagonal matrix of TE modal admittances, and \(
\mathbf{Z}_{s}\in\mathbb{C}^{(2N+1)\times (2N+1)}
\) is the Toeplitz surface-impedance matrix constructed from the Fourier coefficients $z_p$.

Once $\boldsymbol{\Gamma}$ is computed, the Floquet reflection coefficient $B_n$ follow directly from \eqref{eq:bGamma}.
These coefficients characterize the redistribution of the incident power among the reflected Floquet harmonics and form the basis for evaluating the reflection behavior of the RIS.

To characterize the angular distribution of the field reradiated by the RIS, we evaluate its response in the far field. In this region, the distance to an observation point $\textbf{r}$ satisfies $|\mathbf{r}| \gg \lambda$, $|\mathbf{r}| \gg L$, and $L^2/|\mathbf{r}| \ll \lambda$
where $L=\max(2L_x,2L_y)$ denotes the largest size of the RIS \cite{diaz2021macroscopic}. Under these far‑field conditions, the reflected electric field can be written as 
\begin{equation}
E_{\mathrm{ref}}(\theta)
=
\frac{jk}{4\pi}
\frac{e^{-jk|\mathbf{r}|}}{|\mathbf{r}|}
\, A_0 S\, F(\theta)
\label{eq:Escz_compact_ieee}
\end{equation}
where $\theta$ is the spherical-coordinate polar angle of the observation point $\mathbf{r}$, $S=4L_x L_y$ is the area of the RIS, and $F(\theta)$ collects the angular contributions of all reflected Floquet harmonics given by \cite{diaz2021macroscopic} 
\begin{align}
F(\theta)
&=
(\cos\theta - \cos\theta_{\mathrm{i}})
\, \mathrm{sinc}(k \ L_{y \mathrm{ef}})
\notag\\
&\quad+
\sum_{n}
B_n
(\cos\theta + \cos\theta_{\mathrm{r},n})
\, \mathrm{sinc}(k \ L_{y \mathrm{ef},n}),
\label{eq:Ftheta_ieee_final}
\end{align}
where $L_{y\mathrm{ef}}
= L_y(\sin\theta - \sin\theta_{\mathrm{i}})$,
$L_{y\mathrm{ef},n}
= L_y(\sin\theta - \sin\theta_{\mathrm{r},n})$ define the effective size of the RIS for each reflected propagating mode, and $\mathrm{sinc}(x)$ is the $\mathrm{sinc}$ function ($\mathrm{sinc}(x)=\sin(x)/x$). Using the far-field Poynting-expression in \cite{balanis2016antenna}, the radiated power per unit solid angle is 
\begin{equation}
P_{\mathrm{rad}}(\theta)
=
\frac{|\mathbf{r}|^{2}}{2\eta_{0}}\,
\big|E_{\mathrm{ref}}(\theta)\big|^{2}
\label{eq:Prad_E2}
\end{equation}
Substituting~\eqref{eq:Escz_compact_ieee} into~\eqref{eq:Prad_E2} gives 
\begin{equation}
P_{\mathrm{rad}}(\theta)
=
\frac{k^2 |A_0|^2 S^2}{32\pi^2\eta_0}
\, |F(\theta)|^2
\label{eq:Prad_ieee_final}
\end{equation}

Finally, the normalized radiation pattern is defined as
\begin{equation}
\tilde{P}_{\mathrm{rad}}(\theta)
=
\frac{|F(\theta)|^2}
{\displaystyle \max_{\varphi}\, |F(\varphi)|^2}
\label{eq:Pnorm_ieee_final}
\end{equation}
which captures the distribution of power among the reflected modes and provides a direct measure of how the RIS controls the directions of the reflected power.


\section{Numerical Results}

In this section, we evaluate the far‑field behavior of an RIS using the electromagnetic model presented in the previous section. The simulation parameters are listed in Table~\ref{tab:sim_params}.
The RIS operates at $28$~GHz and is illuminated by a normally incident plane wave ($\theta_i=0^\circ$). The surface impedance is designed to steer the reflected wave towards $\theta_r=70^\circ$, and the surface impedance is sampled over one period using $2N{+}1$ points.
\begin{table}[H]
\centering
\caption{Simulation Parameters for the Far field radiation pattern}
\label{table2}
\begin{tabular}{l c}
\toprule
\textbf{Parameter} & \textbf{Value} \\
\midrule
Operating frequency $f$ & $28\,\mathrm{GHz}$ \\
Wavelength & $\lambda = c/f$ \\
Incident angle & $\theta_{\mathrm{i}} = 0^\circ$ \\
Reflection angle & $\theta_{\mathrm{r}} = 70^\circ$ \\
Truncation order & $N = 30$ \\
RIS size & $5D \times 5D$\\
Far-field angle range & $\theta \in [-90^\circ,90^\circ]$\\
Impedance profiles & $Z_1$, $Z_2$, $Z_3$ \\
\bottomrule \label{tab:sim_params}
\end{tabular}
\end{table}
Three commonly used impedance profiles are considered, which are obtained from different design principles and lead to different implementation complexities \cite{di2022communication}. The cotangent impedance defined as
\begin{equation}
    Z_1(y) = j Z_0 \cot\!\big( \tfrac{k}{2}(\sin\theta_i - \sin\theta_r)\, y \big)
\end{equation}
where $Z_{0}=\eta_{0}/\cos\theta_{i}$. This impedance is purely reactive. The geometric-optics (GO) impedance defined as
\begin{equation}
    Z_{2}(y)
= \eta_{0}\,
\frac{1+\Psi(y)}
{\cos\theta_{i}-\Psi(y)\cos\theta_{r}}
\end{equation}
where $\Psi(y)
= e^{-j k (\sin\theta_{r}-\sin\theta_{i})\, y }$. This impedance introduces some power losses, as its real part is typically positive. The globally-optimal impedance
\begin{equation}
    Z_{3}(y)
=\frac{\eta_{0}}{\sqrt{\cos\theta_{i}\cos\theta_{r}}}
\frac{\sqrt{\cos\theta_{r}}+\sqrt{\cos\theta_{i}}\,\Psi(y)}{\sqrt{\cos\theta_{i}}
- \sqrt{\cos\theta_{r}}\,\Psi(y)
}
\end{equation}
which requires local power amplification and losses along the RIS, and is, hence, the most complex design to implement.
For each impedance profile, the reflected Floquet amplitudes are computed using the mode-matching formulation presented in Section. III. 
\begin{figure}[h]
    \centering   \includegraphics[width=0.8\linewidth]{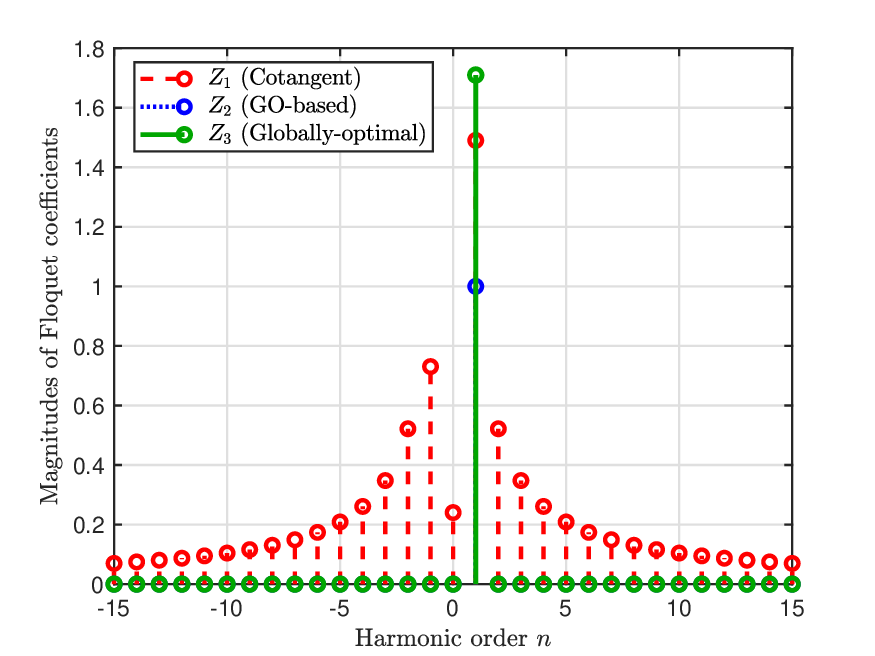}
    \caption{Comparison of the Floquet coefficients magnitude for the three impedance 
profiles.}    \label{fig:harmonic_amp}
\end{figure}

Fig.~\ref{fig:harmonic_amp} compares the magnitudes of the normalized reflected Floquet amplitudes 
$\lvert B_n\rvert$ obtained for the three impedance profiles. For the cotangent impedance $Z_1$, a non-negligible amount of power is distributed across multiple harmonics, indicating that the RIS is unable to effectively concentrate the reradiated power towards the desired direction of reflection that is given by $n=1$. In contrast, both the GO-based impedance
$Z_2$ and the globally-optimal impedance $Z_3$ direct all the reflected power into the desired direction of the reflection, i.e., the harmonic order $n=1$. $Z_3$ gives the strongest contribution while suppressing the power reflected towards any unwanted harmonic orders. This shows that enforcing global power conservation yields a more efficient redistribution of the incident power towards the desired direction.
\vspace{0.7in}

\begin{figure}[h]
    \centering    \includegraphics[width=0.88\linewidth]{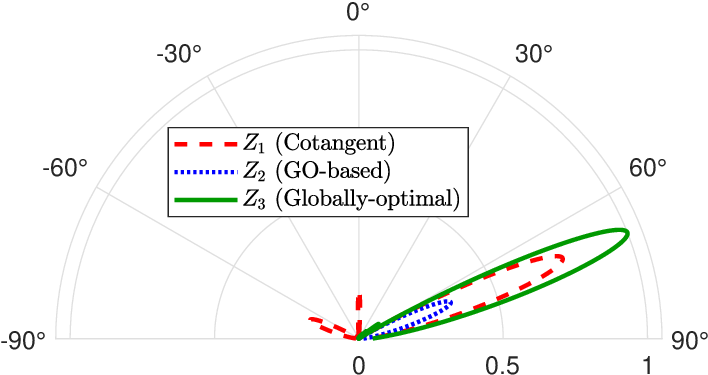}
    \caption{Comparison of the normalized far-field radiation patterns for the three impedance profiles.}
    \label{fig:polar_pattern}
\end{figure}

Figure~\ref{fig:polar_pattern}
compares the normalized far-field radiation patterns for the three impedance profiles. The cotangent impedance $Z_1$ produces a focused beam towards the desired angle $\theta_r=70^\circ$,
but suffers from
the presence of two strong sidelobes. This behavior is consistent with the large number of parasitic Floquet harmonics observed in Fig.~\ref{fig:harmonic_amp}, which translates into undesired radiation at their corresponding angles $\theta_{r,n}$. The GO-based impedance $Z_2$ forms a main lobe with no sidelobes, but radiates less power towards the target direction. The globally-optimal impedance $Z_3$ provides the strongest and the most focused beam at
at the desired direction, with no radiation sidelobes. This behavior matches the harmonic illustration in Fig.~\ref{fig:harmonic_amp}, where the entire reflected power is directed towards the desired mode, i.e., the desired angle $\theta_r$. These results highlight the tradeoff between performance and implementation complexity, and show that the globally-optimal impedance profile, despite its higher complexity, provides better control of the reflected power, by directing it more efficiently towards the desired angle of propagation.

\section{Conclusion}
We presented a physically consistent framework for analyzing anomalous reflection from spatially periodic RIS using the mode-matching technique and Floquet expansion. The results show a clear performance vs. complexity tradeoff where simple impedance profiles are simple to implement but provide limited efficiency, while globally-optimal impedance designs achieve optimal reflection performance at the cost of higher hardware complexity. This tradeoff should guide practical RIS-aided wireless systems design.

\bibliographystyle{IEEEtran}
\bibliography{referencees}

\end{document}